\documentclass{amsart}
\usepackage{amssymb}
\begin{document}
\newtheorem{theorem}{Theorem}[section]
\newtheorem{lemma}[theorem]{Lemma}
\newtheorem{definition}[theorem]{Definition}
\newtheorem{corollary}[theorem]{Corollary}
\newtheorem{remark}[theorem]{Remark}
\def\BB{{\mathcal{B}}}
\def\Tr{\operatorname{Tr}}
\def\BOX{\hbox{$\rlap{$\sqcap$}\sqcup$}}
\def\Pspec{\operatorname{Spec}}
\def\vol{\operatorname{Vol}}
\def\II{\mathcal{I}}
\def\binom#1#2{\!\left(\!\!\!\!\begin{array}{l}#1\\#2\end{array}\!\!\!\!\right)\!}
\makeatletter
  \renewcommand{\theequation}{%
   \thesection.\arabic{equation}}
  \@addtoreset{equation}{section}
 \makeatother
\title[Kaehler Condition]{Spectral Geometry and the Kaehler Condition for Hermitian Manifolds with
Boundary}
\author{JeongHyeong
Park}\thanks{This work was supported by Korea Research
Foundation Grant (KRF-2000-015-DS0003)}
\begin{address}{JHP Dept. of Computer \& Applied Mathematics, Honam
University, Seobongdong 59, Gwangsanku,  Gwangju, 506-714 South
Korea.\quad
Email:jhpark@honam.ac.kr}\end{address}
\begin{abstract} Let $(M,g,J)$ be a compact Hermitian manifold with a smooth boundary. Let
$\Delta_{p,\BB}$ and $\BOX_{p,\BB}$ be the realizations of the real and complex Laplacians on $p$
forms with either Dirichlet or Neumann boundary conditions. We generalize previous results in the
closed setting to show that $(M,g,J)$ is Kaehler if and only if
$\Pspec(\Delta_{p,\BB})=\Pspec(2\phantom{.}\BOX_{p,\BB})$ for $p=0,1$. We also give a characterization of
manifolds with constant sectional curvature or constant Ricci tensor (in the real setting) and manifolds of constant holomorphic
sectional curvature (in the complex setting) in terms of spectral geometry.
\end{abstract}
\subjclass [2000] {58J50}
\keywords{Dirichlet boundary conditions, Neumann boundary conditions, Kaehler, Hermitian manifold,
spectral geometry, constant sectional curvature, constant holomorphic sectional curvature, Einstein manifold}
\maketitle

\section{Introduction}\label{Sect1}

The relationship between the spectrum of certain natural operators of Laplace type
and the underlying geometry of a Riemannian manifold has been studied by many authors. Let
$(M,g)$ be a compact Riemannian manifold with smooth boundary $\partial M$. Let $V$ be a smooth
Hermitian vector bundle over $M$ and let $D$ be a formally self-adjoint operator of Laplace type
acting on the space of smooth sections $C^\infty(V)$. Let $D_\BB$ denote the realization of
$D$ with respect to either the Dirichlet ($\BB=\BB_D)$ or the Neumann ($\BB=\BB_N$) boundary
operators. Then
$D_\BB$ is self-adjoint and has a complete discrete spectral resolution
$\mathcal{S}(D_\BB)=\{(\phi_\nu,\lambda_\nu)\}$. The
$\phi_\nu\in C^\infty(V)$ form a complete orthonormal basis for$L^2(V)$ such that
$$D\phi_\nu=\lambda_\nu\phi_\nu\quad\text{and}\quad\BB\phi_\nu=0.$$
We let the spectrum $\Pspec(D_\BB)=\{\lambda_\nu\}$ be the collection of eigenvalues. We repeat the
eigenvalues according to multiplicity and order the eigenvalues so $\lambda_1\le\lambda_2....$. For
example, if
$M=[0,\pi]$ and if $D=-\partial_x^2$ on $C^\infty(M)$, then:
\begin{eqnarray*}
&&\mathcal{S}(D_{\BB_D})=\left\{\left(\textstyle\sqrt{\frac2\pi}\sin(nx),
    n^2\right)\right\}_{n=1}^\infty\\
&&\Pspec(D_{\BB_D})=\{1,4,9,...\}\\
&&\mathcal{S}(D_{\BB_N})=\left\{\left(1,\textstyle\sqrt{\frac1\pi}\right)\right\}\cup
    \left\{\left(\textstyle\sqrt{\frac2\pi}\cos(nx),n^2\right)\right\}_{n=1}^\infty\\
&&\Pspec(D_{\BB_N})=\{0,1,4,9,...\}.\end{eqnarray*}

Let $(M,g,J)$ be a Hermitian manifold of complex dimension $\hat m$ and corresponding real dimension
$m=2\hat m$; here
$J$ is an integrable almost complex structure which is unitary with respect to the Riemannian metric
$g$. Let
$\Lambda^nM$ be the bundle of complex $n$ forms on $M$. Let
$$\Delta_n:=d^*d+dd^*\quad\text{and}\quad
\BOX_n:=\bar\partial^*\bar\partial+\bar\partial\bar\partial^*\quad\text{on}\quad
C^\infty(\Lambda^nM)$$
be real and complex form valued Laplacians. We further decompose
$$\BOX_n=\oplus_{p+q=n}\phantom{.}\BOX_{(p,q)}\quad\text{on}\quad C^\infty(\Lambda^{(p,q)}M).$$

We introduce the associated {\it Kaehler form} $\Omega(X,Y):=g(X,JY)$. Extend the metric $g$ to be Hermitian on the complexified tangent
bundle. Let $\nabla$ be the Levi-Civita connection.  The following notions are
equivalent and any defines the notion of a {\it Kaehler manifold}:
\begin{enumerate}
\item For every $P$ in $M$, there exist local holomorphic coordinates so
$dg(P)=0$.
\item We have $d\Omega=0$.\item We have $\nabla J=0$.\end{enumerate}Let $\mathcal{I}:=\text{int}(\Omega)$. Let $\delta^\prime$ be the formal adjoint of $\partial$ and$\delta^{\prime\prime}$ be the formal adjoint of $\bar\partial$. For a Kaehler manifold, one has thefollowing relationships:\begin{enumerate}\item $\bar\partial\II-\II\bar\partial=\sqrt{-1}\phantom{.}\delta^\prime$.\item $\bar\partial\delta^\prime+\delta^\prime\bar\partial=0$ and$\partial\delta^{\prime\prime}+\delta^{\prime\prime}\partial=0$.\item $d^*d+dd^*=\partial\delta^\prime+\delta^\prime\partial+\bar\partial\delta^{\prime\prime}+\delta^{\prime\prime}\bar\partial$.\item$\partial\bar\partial\II-\partial\II\bar\partial=\sqrt{-1}\phantom{.}\partial\delta^\prime$and $\bar\partial\II\partial-\II\bar\partial\partial=\sqrt{-1}\phantom{.}\delta^\prime\partial$.
\item
$\sqrt{-1}\phantom{.}(\partial\delta^\prime+\delta^\prime\partial)=\partial\bar\partial\II-\partial\II\bar\partial
+\bar\partial\II\partial-\II\bar\partial\partial$.
\item
$\partial\delta^\prime+\delta^\prime\partial
=\bar\partial\delta^{\prime\prime}+\delta^{\prime\prime}\bar\partial$.
\end{enumerate}

The following well known result is now immediate.

\begin{theorem}\label{thm1.1} Let $(M,g,J)$ be a compact Kaehler manifold without boundary of complex
dimension
$\hat m$. Then $\Delta=2\phantom{.}\BOX$ and
so $\Pspec(\Delta_p)=\Pspec(2\phantom{.}\BOX_p)$ for all $p$.\end{theorem}

Conversely, one has the following result to T. Tsujishita (reported by Gilkey \cite{Gi75}):

\begin{theorem}\label{thm1.2} Let $(M,g,J)$ be a compact Hermitian manifold without
boundary. If $\Pspec(\Delta_0)=\Pspec(2\BOX_{0})$ and if $\Pspec(\Delta_1)=\Pspec(2\BOX_1)$, then
$(M,g,J)$ is Kaehler.
\end{theorem}

Donnelly \cite{Do75} established a similar characterization of the Kaehler property using the reduced
complex Laplacian. Pak
\cite{refHKPak02} extended these results to the context of almost isospectral manifolds.\medbreak

 Theorem \ref{thm1.2} is sharp.  We refer to Gilkey \cite{Gi74} for the
proof of the following result:

\begin{theorem}\label{thm1.3} Let
$ds^2=dz_1\circ d\bar z_1+e^{\psi(z_1)}dz_2\circ d\bar z_2+e^{-\psi(z_1)}dz_3\circ d\bar z_3$
be a Hermitian metric on the torus $M_3$ where $\psi(z_1)$ is an arbitrary smooth real valued
function. Then $\Delta_0=2\phantom{.}\BOX_{0}$ but the metric is not Kaehler.\end{theorem}

A Riemannian manifold of constant sectional curvature $c$ is said to be a {\it space form}; a
Kaehler manifold of constant holomorphic sectional curvature $c$ is said to be a {\it
complex space form}. Modulo rescaling, any space form is locally isometric to the unit sphere, to
flat space, or to hyperbolic space. Similarly, modulo rescaling, any complex space form is locally
isometric to complex projective space, to flat space, or to the negative curvature dual. Thus the
geometries are very rigid in this context.

Patodi \cite{refPa70} established the following spectral characterization of space forms:
\begin{theorem}\label{thm1.4}   Let $(M_i,g_i)$ be compact Riemannian manifolds
without boundary. Assume that $\Pspec(\Delta_{p},M_1)=\Pspec(\Delta_p,M_2)$ for $0\le p\le 2$.
Then:\begin{enumerate}\item The manifold $M_1$ has constant scalar curvature $c$ if and only if the manifold $M_2$ has constant scalar curvature $c$.
\item The manifold $M_1$ is Einstein if and only if the manifold $M_2$ is Einstein.
\item The manifold $M_1$ has constant sectional curvature $c$ if and only if the manifold $M_2$ has
constant sectional curvature $c$.
\end{enumerate}\end{theorem}

 Donnelly \cite{Do75a} and Gilkey and Sacks
\cite{GiSa75} extended Theorem \ref{thm1.4} to the complex
setting -- see also related work by Friedland \cite{refFri95,refFri96}, C.C.
Hsuing et. al. \cite{Hsuing96}, and Pak \cite{refHKPak03}.

\begin{theorem}\label{thm1.5}  Let $(M_i,g_i,J_i)$ be compact Kaehler manifolds
without boundary. Assume that
$\Pspec(\BOX_{p,q},M_1)=\Pspec(\BOX_{p,q},M_2)$ for $0\le p\le 2$ and $0\le q\le 2$. Then the manifold
$M_1$ has constant holomorphic sectional curvature $c$ if and only if the manifold $M_2$ has constant
holomorphic sectional curvature $c$.
\end{theorem}

\medbreak We can  extend Theorem \ref{thm1.2} to  manifolds
with boundary:

\begin{theorem}\label{thm1.6} Let $\BB$ denote either Dirichlet or Neumann boundary conditions. Let
$(M,g,J)$ be a compact Hermitian manifold with smooth boundary $\partial M$. Assume that
$\Pspec(\Delta_{0,\BB})=\Pspec(2\phantom{.}\BOX_{0,\BB})$ and that
$\Pspec(\Delta_{1,\BB})=\Pspec(2\phantom{.}\BOX_{1,\BB})$.
Then
$(M,g,J)$ is Kaehler.
\end{theorem}

We can also generalize Theorems \ref{thm1.4} and \ref{thm1.5} to the category of manifolds with boundary
under the additional technical hypothesis that the manifolds in question have constant scalar curvature.

\begin{theorem}\label{thm1.7} Let $\BB$ denote either Dirichlet or Neumann boundary conditions.
Let $(M_i,g_i)$ be compact Riemannian manifolds
with smooth boundaries $\partial M_i$ which have constant scalar curvatures $\tau_i$. Assume that
$\Pspec(\Delta_{p,\BB},M_1)=\Pspec(\Delta_{p,\BB},M_2)$ for
$0\le p\le 2$. Then:
\begin{enumerate}
\item $\tau_1=\tau_2$.
\item The manifold $M_1$ is Einstein if and only if the manifold $M_2$ is Einstein.
\item The manifold $M_1$ has constant sectional curvature $c$ if and only if the manifold $M_2$ has
constant sectional curvature $c$.
\end{enumerate}\end{theorem}

\begin{theorem}\label{thm1.8}Let $\BB$ denote either Dirichlet or Neumann boundary conditions.
 Let $(M_i,g_i,J_i)$ be compact Kaehler manifolds with smooth
boundaries $\partial M_i$ which have constant scalar curvatures $\tau_i$. Assume that
$\Pspec(\BOX_{p,q},M_1)=\Pspec(\BOX_{p,q},M_2)$ for
$0\le p\le 2,0\le q\le2$. Then the manifold $M_1$ has
constant holomorphic sectional curvature $c$ if and only if the manifold $M_2$ has constant
holomorphic sectional curvature $c$.
\end{theorem}

Here is a brief outline to the remainder of this paper. In Section \ref{Sect2}, we review some previous
results concerning the heat trace asymptotics. In Section \ref{Sect3}, we complete the
proof of Theorem \ref{thm1.6} and in Section \ref{Sect4}, we complete the proof of Theorems
\ref{thm1.7} and
\ref{thm1.8}. In Remark \ref{rmk-3}, we extend Theorems \ref{thm1.7} and \ref{thm1.8} to more general boundary
conditions of Robin type where the auxiliary endomorphism $S$ is `universal' in a certain sense.

It is a pleasant task to thank the referee for helpful suggestions concerning the manuscript.

\section{Heat trace asymptotics}\label{Sect2}

Let $M$ be a compact Riemannian manifold of real dimension $m$ with smooth boundary $\partial M$ and
let $D_\BB$ be the realization of an operator of Laplace type on $M$ with respect to either Dirichlet
or Neumann boundary conditions. Let $e^{-tD_\BB}$ be the fundamental solution of the heat equation.
This operator is of trace class and as $t\downarrow0$ there is a complete asymptotic expansion with
locally computable coefficients in the form:
$$\Tr _{L^2}e^{-tD_\BB}\sim\textstyle\sum_{n\ge0}t^{(n-m)/2}a_n(D,\BB).$$

To study the heat trace coefficients $a_n(D,\BB)$, we must introduce a bit of additional notation. There
is a canonically defined connection
$\nabla=\nabla(D)$ and a canonically defined endomorphism $E=E(D)$ so that
$$D=-(\Tr(\nabla^2)+E).$$

Let indices $i$, $j$, $k$ range from $1$ to $m$ and index a local orthonormal frame $\{e_1,...,e_m\}$
for
$TM$. Let indices $a$, $b$, and $c$ range from $1$ to $m-1$ and index a local orthonormal frame
$\{e_1,...,e_{m-1}\}$ for
$T\partial M$; on $\partial M$, we let $e_m$ be the inward unit normal vector field. Let
$\Omega$ be the curvature of
$\nabla$, let
$\tau:=R_{ijji}$ be the normalized scalar curvature, let $\rho_{ij}:=R_{ikkj}$ be the Ricci tensor, and
let $L_{ab}$ be the second fundamental form. We adopt the Einstein convention and sum over repeated
indices. Let `;' denote multiple covariant differentiation. We refer to
\cite{BrGi} for the proof of the following result:
\begin{theorem}\label{thm2.1} Let $D$ be an operator of Laplace type on the space of sections
$C^\infty(V)$ to a vector bundle $V$ over a compact manifold
$M$ with smooth boundary $\partial M$. Let $I$ be the identity endomorphism of $V$. With Dirichlet boundary conditions, we have:
\begin{enumerate}
\item $a_{0}(D,\BB_D)=(4\pi)^{-m/2}\int_M\Tr \{I\}$.
\smallbreak\item
$a_{1}(D,\BB_D)=-(4\pi)^{-(m-1)/2}\frac14\int_{\partial M}\Tr \{I\}$.
\smallbreak\item $ a_{2}(D,\BB_D)=(4\pi)^{-m/2}\frac16\int_M\Tr\{6E+\tau I\}$
$+(4\pi)^{-m/2}\frac16\int_{\partial M}\Tr\{2L_{aa} I\}$.
\smallbreak\item $a_{3}(D,\BB_D)=-\frac1{384}(4\pi)^{-(m-1)/2}
\int_{\partial M}\Tr\{96E+(16\tau +8R_{am
am}$\newline$
+7L_{aa}L_{bb}-10L_{ab}L_{ab})I\}$.
\smallbreak\item
$a_{4}(D,\BB_D)=(4\pi)^{-m/2}\frac1{360}\int_M\Tr\{60E_{;kk}+60\tau
E+180E^{2}+30\Omega^{2}$\newline$
+(12\tau_{;kk}+5\tau^{2}-2|\rho|^{2}+2|R|^{2})I\}$
$+(4\pi)^{-m/2}\frac1{360}\int_{\partial M}\Tr\{-120E_{;m}$\newline$+120
EL_{aa}+(-18\tau_{;m}+20\tau L_{aa}$
 $+4R_{am
am}L_{bb} -12R_{am bm}L_{ab}$\newline$+4R_{ab
cb}L_{ac}+\frac{40}{21}L_{aa}L_{bb}L_{cc}
-\frac{88}7L_{ab}L_{ab}L_{cc}+\frac{320}{21}L_{ab}L_{bc}L_{ac})I\}$.
\end{enumerate}
\end{theorem}

Let $\nabla_m$ denote the covariant derivative with respect to $e_m$ on $\partial M$. Let $S$ be an auxiliary
endomorphism of $V|_{\partial M}$. The Robin boundary operator is then given by:
$$\BB_S\phi:=(\nabla_m\phi+S\phi)|_{\partial M}\,.$$
We take $S=0$ to define Neumann boundary conditions. Again, we refer to \cite{BrGi} for the proof of the following result:

\begin{theorem}\label{thm2.2} With Robin boundary conditions, we have:
\begin{enumerate}
\item
$a_0(D,\BB_S)=(4\pi)^{-m/2}\int_{M}\Tr \{I\}$. \smallbreak\item
$a_1(D,\BB_S)=
     (4\pi)^{(1-m)/2}\frac14\int_{\partial M}\Tr \{I\}$.
\smallbreak\item $a_2(D,\BB_S)=(4\pi)^{-m/2}\frac16
     \int_{M}\Tr\{6E+\tau I\}$\newline
$+(4\pi)^{-m/2}\frac16\int_{\partial M}\Tr\{2L_{aa}I+12S\}$.
\smallbreak\item $a_3(D,\BB_S)= (4\pi)^{(1-m)/2}\frac1{384}
\int_{\partial M} \Tr\{96E+(16\tau -8R_{amma}$
$+13L_{aa}L_{bb}$\newline$ +2L_{ab}L_{ab})I+96SL_{aa}+192S^2\}$.
\smallbreak\item $a_4 (D,\BB_S)
=(4\pi)^{-m/2}\frac1{360}\int_{M}\Tr\{60E_{;kk}+60 \tau E
+180E^{2}+30\Omega^2$\newline $ +(12 R_{ijji;kk}$ $ + 5\tau^2 -
2|\rho|^2+2|R|^2)I\}$ $+(4\pi)^{-m/2}\frac1{360}\int_{\partial M}
\Tr\{
      240 E _{;m}$\newline$
         +120 E L_{aa}+ (42 R _{ijji;m}$
$  + 20 \tau L_{aa}
 + 4 R _{am a m } L_{bb} -12 R_{ambm} L_{ab}$\newline$
 + 4 R_{abcb} L_{ac} $
$  + \frac{40}3 L_{aa} L_{bb} L_{cc}
 +8
   L_{ab} L_{ab}  L_{cc}$$ + \frac{32}3 L_{ab}
L_{bc}  L_{ac})I+120S\tau$
\newline$+720SE+144SL_{aa}L_{bb}+48SL_{ab}L_{ab}+480S^2L_{aa}+480S^3\} $.
\end{enumerate}
\end{theorem}

\section{The proof of Theorem \ref{thm1.6}}\label{Sect3}

Let $\star$ be the Hodge operator. We introduce the following
invariants:\begin{enumerate}
\item Let $K_1:=\star(d\bar\partial\Omega\wedge\Omega^{\hat m-2})$ for $\hat m\ge2$.
\item Let $K_2:=\frac12|d\Omega|^2$ for $\hat m\ge2$.
\item Let $K_3:=\star(d\Omega\wedge\bar\partial\Omega\wedge\Omega^{\hat m-3})$ for $\hat m\ge3$.
\item Let $\kappa:=L_{aa}$ be the geodesic curvature of the boundary.
\end{enumerate}

We may then use Theorem \ref{thm2.1} to extend results of Gilkey
\cite{Gi74} to see:
\medbreak\quad
$\textstyle a_2(2\phantom{.}\BOX_{0})=(4\pi)^{-\hat m}\frac16\{\int_{\partial M}2\kappa+
\textstyle\int_M(\tau+3K_2+3K_3)\}$
\medbreak\quad$
\textstyle a_2(\Delta_0)=(4\pi)^{-\hat m}\frac16\{\int_{\partial M}2\kappa+
\textstyle\int_M\tau\}$
\medbreak\quad$
\textstyle a_2(2\phantom{.}\BOX_{1,0})=(4\pi)^{-\hat m}\frac16\{\int_{\partial M}2\hat
m\kappa+ (\hat m-3)\textstyle\int_M(\tau+3K_2+3K_3)$\smallbreak\qquad$+
\textstyle\int_M(-6K_1+6K_2+3K_3)\}$\medbreak\quad$
\textstyle a_2(2\phantom{.}\BOX_{0,1})=(4\pi)^{-\hat m}\frac16\{\int_{\partial M}2\hat
m\kappa+ (\hat m-3)\textstyle\int_M(\tau+3K_2+3K_3)$\smallbreak\qquad$+
\textstyle\int_M(6K_1+6K_2+3K_3)\}$\medbreak\quad$
\textstyle a_2(\Delta_1)=(4\pi)^{-\hat m}\frac16\{\int_{\partial
M}4\hat m\kappa+\textstyle2(\hat m-3)\int_M\tau.$\medbreak\noindent
Since $a_2(\Delta_0)=a_2(2\phantom{.}\BOX_{0})$, we have
\begin{equation}
\textstyle\int_M(3K_2+3K_3)=0.\label{eqn3.4}\end{equation}
We use this relation and the relation $a_2(\Delta_1)=a_2(2\phantom{.}\BOX_1)$ to see
\begin{equation}
\textstyle\int_M(6K_2+3K_3)=0.\label{eqn3.5}\end{equation}
Equations (\ref{eqn3.4}) and (\ref{eqn3.5}) then imply
$\int_MK_2=0$ and hence $M$ is Kaehler.

If $\hat m=2$, we have the formulae:
\begin{equation}\begin{array}{lll}
\textstyle a_2(2\phantom{.}\BOX_{0})&=(4\pi)^{-\hat m}\frac16\{\int_{\partial M}2\kappa&+
\textstyle\int_M(\tau+3K_2)\}\\
\textstyle a_2(\Delta_0)&=(4\pi)^{-\hat m}\frac16\{\int_{\partial M}2\kappa&+
\textstyle\int_M\tau\}.
\end{array}\end{equation}
Since $a_2(\Delta_0)=a_2(2\phantom{.}\BOX_{0})$, $\int_MK_2=0$ and $M$ is Kaehler; the condition
relating $\Delta_1$ and $2\phantom{.}\BOX_1$ is not necessary in this instance. Finally, if $\hat m=1$,
then $M$ is automatically Kaehler. \hfill $\BOX$

\section{Proof of Theorems \ref{thm1.7} and \ref{thm1.8}}\label{Sect4} We set $S=0$.
Let $\varepsilon:=\frac14$ for
Neumann boundary conditions and $\varepsilon:=-\frac14$ for
Dirichlet boundary conditions. By Theorems \ref{thm2.1} and
\ref{thm2.2},
\begin{eqnarray}
&&\textstyle\Tr _{L^2}(e^{-t\Delta_{0,\BB}})=(4\pi
t)^{-m/2}\{\vol(M)+\varepsilon{\sqrt{t}}\vol(\partial M)+O(t)\}\quad\text{so}\nonumber\\
&&\vol(M_1)=\vol(M_2),\quad\vol(\partial M_1)=\vol(\partial M_2),\quad\text{and}\label{eqn4.1}\\
&&\dim_{\mathbb{R}}(M_1)=\dim_{\mathbb{R}}(M_2).
\nonumber\end{eqnarray}
We set $m:=\dim_{\mathbb{R}}(M_i)$ to this common value and compute:
\begin{eqnarray*}
&&a_2(\Delta_{0,\BB},M_i)=(4\pi)^{-m/2}\textstyle\frac16\{\int_{M_i}\tau_i+\int_{\partial
M_i}2\kappa_i\}\\
&&a_2(\Delta_{1,\BB},M_i)=(4\pi)^{-m/2}\textstyle\frac16\{\int_{M_i}(m-6)\tau_i+\int_{\partial
M_i}2m\kappa_i\}.
\end{eqnarray*}
We may then establish assertion (1) by computing:
\begin{eqnarray*}
\textstyle\tau_1&=&(4\pi)^{m/2}\vol(M_1)^{-1}\{ma_2(\Delta_{0,\BB},M_1)-a_2(\Delta_{1,\BB},M_1)\}\\
&=&(4\pi)^{m/2}\vol(M_2)^{-1}\{ma_2(\Delta_{0,\BB},M_2)-a_2(\Delta_{1,\BB},M_2)\}\\
&=&\tau_2\,.\end{eqnarray*}
For subsequent use, we compute similarly that:
\begin{eqnarray}
\textstyle\int_{\partial M_1}\kappa_1&=&
3(4\pi)^{m/2}a_2(\Delta_{0,\BB},M_1)
 -\textstyle\frac12\int_{M_1}\tau_1\nonumber\\
&=&
3(4\pi)^{m/2}a_2(\Delta_{0,\BB},M_2)
 -\textstyle\frac12\int_{M_2}\tau_2\label{eqn4.2}\\
&=&\textstyle\int_{\partial M_2}\kappa_2\,.\nonumber
\end{eqnarray}
\medbreak The interior integrands defining $a_4(\Delta_p)$ have been determined by Patodi \cite{refPa70}. Motivated by his work,
we introduce constants:
\begin{eqnarray*}
&&c_{m,p}^1=\phantom{-.}\textstyle\frac1{72}\frac{m!}{p!(m-p)!}-\phantom{.}\frac16\frac{(m-2)!}{(p-1)!(m-p-1)!}
           +\frac12\frac{(m-4)!}{(p-2)!(m-p-2)!},\\
&&c_{m,p}^2=-\textstyle\frac1{180}\frac{m!}{p!(m-p)!}+\phantom{.}\frac12\frac{(m-2)!}{(p-1)!(m-p-1)!}-2\frac{(m-4)!}{(p-2)!(m-p-2)!},\\
&&c_{m,p}^3=\phantom{-}\textstyle\frac1{180}\frac{m!}{p!(m-p)!}-\frac1{12}\frac{(m-2)!}{(p-1)!(m-p-1)!}
         +\frac12\frac{(m-4)!}{(p-2)!(m-p-2)!},\\
&&c_{m,p}^4=\phantom{-.}\textstyle\frac1{30}\frac{m!}{p!(m-p)!}-\phantom{.}\frac16\frac{(m-2)!}{(p-1)!(m-p-1)!}\,.
\end{eqnarray*}
The work of Patodi then shows if $\partial M$ is empty that:
\begin{equation}\label{jhp1}
a_4(\Delta_p)
  =(4\pi)^{-m/2}\textstyle\int_M\{c_{m,p}^1\tau^2+c_{m,p}^2|\rho|^2+c_{m,p}^3|R|^2+c_{m,p}^4\tau_{;ii}\}\,.
\end{equation}

To simplify the notation, we introduce reduced invariants
$$\tilde a_n(\Delta_p,\BB):=a_n(\Delta_p,\BB)-\textstyle\frac{m!}{p!(m-p)!}a_n(\Delta_0,\BB)\,.$$
The terms
$$\{R_{amma}L_{bb},\ R_{ammb}L_{ab},\ R_{abcb}L_{ac}\ L_{aa}L_{bb}L_{cc},\ L_{ab}L_{ab}L_{cc},\ L_{ab}L_{bc}L_{ac}\}$$
appearing in Theorem \ref{thm2.1} are all multiplied by $\Tr(I_{\Lambda^p})=\textstyle\frac{m!}{p!(m-p)!}$. Thus
they do not appear in
$\tilde a_4(\Delta_{p,\BB})$; only the terms involving $E(\Delta_p)$
survive in the boundary contributions. One can use the Weitzenb\"och formula to see that
$$\Tr(E(\Delta_p))=-\textstyle\frac{(m-2)!}{(p-1)!(m-p-1)!}\tau\,.$$
Using equation (\ref{jhp1}), we see there
exist universal constants so
\begin{eqnarray}
\tilde a_4(\Delta_{p,\BB})\nonumber
&=&\textstyle(4\pi)^{-m/2}\int_{M}
       \{\tilde c_{m,p}^1\tau^2+\tilde c_{m,p}^2|\rho|^2+\tilde c_{m,p}^3|R|^2+
        \tilde c_{m,p}^4\tau_{;kk}\}\label{eqn4x}\\
&+&\textstyle\textstyle(4\pi)^{-m/2}
\int_{\partial M}\{\tilde c_{m,p}^5\tau\kappa+\tilde c_{m,p}^6\tau_{;m}\},\nonumber
\end{eqnarray}
where, by Patodi's result, we have:
\begin{eqnarray}
&&\textstyle\tilde c_{m,p}^1:=c_{m,p}^1-\textstyle\frac{m!}{p!(m-p)!}c_{m,0}^1=-\phantom{.}\frac16\frac{(m-2)!}{(p-1)!(m-p-1)!}
           +\frac12\frac{(m-4)!}{(p-2)!(m-p-2)!},\label{jhp-2}\\
&&\tilde c_{m,p}^2:=c_{m,p}^2-\textstyle\frac{m!}{p!(m-p)!}c_{m,0}^2=
         \phantom{..a}\frac12\frac{(m-2)!}{(p-1)!(m-p-1)!}
         -2\frac{(m-4)!}{(p-2)!(m-p-2)!},\nonumber\\
&&\tilde c_{m,p}^3:=c_{m,p}^3-\textstyle\frac{m!}{p!(m-p)!}c_{m,0}^3=
         -\frac1{12}\frac{(m-2)!}{(p-1)!(m-p-1)!}
         +\frac12\frac{(m-4)!}{(p-2)!(m-p-2)!},\nonumber\\
&&\tilde c_{m,p}^4:=c_{m,p}^4-\textstyle\frac{m!}{p!(m-p)!}c_{m,0}^4
  =-\phantom{.}\frac16\frac{(m-2)!}{(p-1)!(m-p-1)!}\,.\nonumber
\end{eqnarray}

For $p=1,2$, we have, by assumption, that:
\begin{eqnarray}
 &&\tilde a_4(\Delta_{p,\BB},M_1)=a_4(\Delta_{p,\BB},M_1)-\textstyle\frac{m!}{p!(m-p)!}a_4(\Delta_{0,\BB},M_1)\label{eqn4xy}\\
&=&\tilde a_4(\Delta_{p,\BB},M_2)=a_4(\Delta_{p,\BB},M_2)-\textstyle\frac{m!}{p!(m-p)!}a_4(\Delta_{0,\BB},M_2)\,.\nonumber
\end{eqnarray}
By equation (\ref{eqn4.2}), $\int_{\partial M_1}\kappa_1=\int_{\partial M_2}\kappa_2$. By assertion (1),
$\tau_1=\tau_2$. Since the scalar curvature is constant, $\tau_{;m}=0$. Thus since $\vol(\partial M_1)=\vol(\partial
M_2)$,  the boundary integrals are equal. Furthermore, since
$\vol(M_1)=\vol(M_2)$, the interior integrals of
$\tau^2$ are equal. Since $\tau_{;ii}=0$, we have
\begin{equation}\label{jhp-eqn4}
\textstyle\int_{M_1}(\tilde c_{m,p}^2|\rho_1|^2+\tilde c_{m,p}^3|R_1|^2)
=\textstyle\int_{M_2}(\tilde c_{m,p}^2|\rho_2|^2+\tilde c_{m,p}^3|R_2|^2)
\end{equation}
for $n=1,2$; these two equations are independent since, by display (\ref{jhp-2}),
$$\det\left(\begin{array}{cc}\tilde c_{m,1}^2&\tilde c_{m,1}^3\\
\tilde c_{m,2}^2&\tilde c_{m,2}^3\vphantom{\vrule height 14pt}\end{array}\right)=
\det\left(\begin{array}{cc}\textstyle\frac12&-\frac1{12}\\\textstyle\frac{m-2}2-2&-\frac{m-2}{12}+\frac12
\vphantom{\vrule height 14pt}
\end{array}\right)=\textstyle\frac14-\frac16\ne0\,.$$
Consequently
\begin{equation}\textstyle\int_{M_1}|\rho_1|^2=\int_{M_2}|\rho_2|^2\quad\text{and}\quad
\int_{M_1}|R_1|^2=\int_{M_2}|R_2|^2.\label{eqn443}\end{equation}

A manifold $M$ has constant sectional curvature $c$ if and only if
$$
0=\textstyle\int_M|R_{ijkl}-c(\delta_{il}\delta_{jk}-\delta_{ik}\delta_{jl})|^2\\
=\textstyle\int_M(|R|^2-4c\tau+c^2\varepsilon_m)$$
where $\varepsilon_m:=|\delta_{il}\delta_{jk}-\delta_{ik}\delta_{jl}|^2$ is polynomial in $m$.
We use equation (\ref{eqn4.1}), equation (\ref{eqn443}), and assertion (1) to complete the proof of Theorem
\ref{thm1.7} (3) by computing:
$$\textstyle\int_{M_1}(|R_1|^2-4c\tau_1+c^2\varepsilon_m)=
\textstyle\int_{M_2}(|R_2|^2-4c\tau_2+c^2\varepsilon_m)\,.$$
Note that $M$ is Einstein if and only if there is a constant $c$ so
$$0=|\rho_{ij}-c\delta_{ij}|^2=|\rho|^2-2c\tau+mc^2\,.$$
Thus Theorem \ref{thm1.7} (2) can be established by verifying that:
\begin{eqnarray*}
0&=&\textstyle\int_{M_1}(|\rho_1|^2-2c\tau_1+mc^2)=\textstyle\int_{M_2}(|\rho_2|^2-2c\tau_2+mc^2)\,.
\end{eqnarray*}

\medbreak As a similar argument based on the results of
\cite{Do75a,GiSa75} establishes Theorem \ref{thm1.8}, we shall omit the details of the proof of Theorem
\ref{thm1.8} in the interests of brevity.

\begin{remark}\label{rmk-3}\rm We can generalize Theorems \ref{thm1.7} and \ref{thm1.8} to the context of Robin
boundary conditions as follows. One could take $S=c_1+c_2\kappa$; the same cancellation
argument as that given above to establish equation (\ref{jhp-eqn4}) shows the additional boundary terms cancel off for the reduced
invariant. What is crucial is that the boundary condition be natural and universal in the context in which we are working.
\end{remark}

\end{document}